\documentclass[11pt]{article}

\usepackage{latexsym,mathrsfs}
\usepackage{amsthm,enumerate,verbatim}
\usepackage{amsfonts}
\usepackage{url}
\usepackage[utf8]{inputenc}

\usepackage{hyperref}

\usepackage{amsmath}
\usepackage{cases}
\usepackage{amssymb}
\usepackage{paralist}
\usepackage{multirow}
\usepackage{array}
\usepackage{graphicx}
\usepackage{threeparttable}

\usepackage{algorithm}
\usepackage{algorithmic}

\newtheorem*{definition*}{Definition}

\newtheorem{remark}{Remark}

\DeclareMathOperator{\rank}{rank}

\title{ Improved SVD-based Initialization for Nonnegative Matrix Factorization using Low-Rank Correction } 
\date{}

\author{Atif Muhammad Syed\textsuperscript{1}\,, \, Sameer Qazi\textsuperscript{1}, \, Nicolas Gillis\textsuperscript{2,}\thanks{Corresponding author. Email: nicolas.gillis@umons.ac.be. }\\ \textsuperscript{1}Graduate School of Science and Engineering,\\ PAF-Karachi Institute of Economics and Technology\\Karachi, Pakistan\\
\textsuperscript{2}Department of Mathematics and Operational Research\\Faculté polytechnique, Université de Mons\\ Rue de Houdain 9, 7000 Mons, Belgium}

\begin{document}

\maketitle

\begin{abstract}

Due to the iterative nature of most nonnegative matrix factorization (\textsc{NMF}) algorithms, initialization is a key aspect as it significantly influences both the convergence and the final solution obtained. Many initialization schemes have been proposed for NMF, among which one of the most popular class of methods are based on the singular value decomposition (SVD). 	However, these SVD-based initializations do not satisfy a rather natural condition, namely that the error should decrease as the rank of factorization increases. In this paper, we propose a novel SVD-based \textsc{NMF} initialization to specifically address this shortcoming by taking into account the SVD factors that were discarded to obtain a nonnegative initialization. 	This method, referred to as nonnegative SVD with low-rank correction (NNSVD-LRC), allows us to significantly reduce the initial error at a negligible additional computational cost using the low-rank structure of the discarded SVD factors. NNSVD-LRC has two other advantages compared to previous SVD-based initializations: (1) it provably generates sparse initial factors, and (2) it is faster as it only requires to compute a truncated SVD of rank $\lceil r/2 + 1 \rceil$ where $r$ is the factorization rank of the sought NMF decomposition (as opposed to a rank-$r$ truncated SVD for other methods). We show on several standard dense and sparse data sets that our new method competes favorably with state-of-the-art SVD-based initializations for NMF. 

\end{abstract}


	\section{Introduction}
	
	Nonnegative matrix factorization (NMF) is the problem of approximating a input nonnegative matrix $X$ as the product of two nonnegative matrices: Given $X \in\mathbb R_{\geq0}^{m\times n}$ and an integer~$r$, find $W\in\mathbb R_{\geq0}^{m\times r}$ and $H\in\mathbb R_{\geq0}^{r\times n}$ such that $X\approx WH$. 
	\textsc{NMF} allows to reconstruct data using a purely additive model: each column of $X$ is a nonnegative linear combination of the columns of $W$. 
	For this reason, it is widely employed in research fields like image processing and 
	computer vision~\cite{ensari_character_2016, prajapati_brain_2015}, 
	data mining and document clustering~\cite{du_hybrid_2017}, 
	hyperspectral image analysis~\cite{luce_using_2016, shiga_sparse_2016}, 
	signal processing~\cite{yoshii_students_2016} and 
	computational biology~\cite{maruyama_detecting_2014}; 
	see also~\cite{cichocki2009nonnegative, gillis2014} and the references therein. 
	
	To measure the quality of the \textsc{NMF} approximation, a distance metric should be chosen. In this paper, we focus on the most widely used one, namely the Frobenius norm, leading to the following optimization problem 	
	\begin{equation}
	\label{eq:nmf_min_prob} 
	\min_{ W \in \mathbb R^{m \times r}, H \in \mathbb R^{r \times n} }  
	\| X-WH \| _F^2  \quad \text{such that} \quad W\geq0 \text{ and } H\geq0, 
	\end{equation}
	where $\|M\|_F=\sqrt{\sum_{i,j} M^{2}_{i,j}}$  is Frobenius norm of a matrix $M$. 
	Most algorithms tackling~\eqref{eq:nmf_min_prob} use standard non-linear optimization schemes such as block coordinate descent methods hence initialization of the factors $(W,H)$ is crucial in practice as it will influence 
	\begin{itemize}
	    \item[(i)] the number of iterations needed for an algorithm to converge (in fact, if the initial point is closer to a local minimum, it will require less iterations to converge to it), 
	and
	\item[(ii)] the final solution to which the algorithm will converge. 
	\end{itemize} 
	Many approaches have been proposed for \textsc{NMF} initialization, for example based 
	on $k$-means and spherical $k$-means~\cite{wild2004improving}, 
	on fuzzy $c$-means~\cite{rezaei_efficient_2011}, 
	on nature inspired heuristic algorithms~\cite{tan_using_2011}, 
	on Lanczos bidiagonalization~\cite{wang_effective_2012}, 
	on subtractive clustering~\cite{casalino_subtractive_2014},
	and on the successive projection algorithm~\cite{sauwen_successive_2017}, to name a few; see also~\cite{langville2006initializations}. 
	
In this paper, we focus on SVD-based initializations for NMF. 
			Two of the most widely used methods are \textsc{NNDSVD} \cite{boutsidis_svd_2008} and \textsc{SVD-NMF} \cite{qiao_new_2015} which are described in the next section. 
			These methods suffer from the fact that the approximation error $||X-WH||_F^2$ of the initial factors $(W,H)$ increases as the rank increases which is not a desirable property for NMF initializations.   
		Our key contribution is to provide a new SVD-based initialization that does not suffer from this shortcoming while 
		(i) it generates sparse factors which not only provide storage efficiency~\cite{gillis_using_2010} but also provide better part-based  representations~\cite{casalino_subtractive_2014,elad_role_2010} and resilience to noise~\cite{ye_multitask_2015,sun_graph_2016}, 
		and 
		(ii) it only requires a truncated SVD of rank $\lceil \frac{r}{2} + 1 \rceil$, as opposed to a truncated SVD of rank $r$ for the other SVD-based initializations. 
		
	
	\paragraph{Outline of the paper} This paper is organized as follows. 
	Section~\ref{proposedsol} will discuss our proposed solution in details, 
	highlighting the differences with existing SVD-based initializations. 
	In Section~\ref{numexp}, we evaluate our proposed solution against other SVD-based initializations on dense and sparse data sets. 
	Section~\ref{concl} concludes the paper.

	\section{Nonnegative SVD with low-rank correction, a new SVD-based NMF initialization} \label{proposedsol}
	
		The truncated SVD is a low-rank matrix approximation technique that approximates a given matrix \(X\in \mathbb{R}^{m\times n}\)  as a sum of \(r\) rank-one terms made of singular triplets, where \(1 \leq r \leq \rank(X)\). 
		Each singular triplet $(u_i,v_i,\sigma_i)$ ($1 \leq i \leq r$) consists of two column vectors \(u_i\) and \(v_i\) which are the left and the right singular vectors, respectively, associated with the $i$th singular value (which we assume are sorted in nonincreasing order).  We have 
		\begin{equation}
		X\approx X_{r} = \sum_{i=1}^{r} \sigma_{i}u_{i}v_{i}^{T} = U_{r}\Sigma_{r} V_{r}^{T}, 
		\end{equation}
where 
		$(.)^T$ is the transpose of given matrix or vector, 
		$X_{r}$ is the rank-$r$ approximation of $X$, 
		the columns of 
		$U_{r} \in \mathbb{R}^{m\times r}$ (resp.\@ of $V_{r} \in \mathbb{R}^{n \times r}$) are the left (resp.\@ right) singular vectors,  
		and $\Sigma_{r} \in \mathbb{R}^{r\times r}$ is the diagonal matrix containing the singular values on its diagonal. 
		According to Eckhart-Young theorem, $X_r$ provides an optimal rank-$r$ approximation of \(X\) with respect to the Frobenius and spectral norms~\cite{golub_matrix_2013}. 
		To simplify our later derivations, we transform the three factors of the SVD representation into two factors, like in NMF, by multiplying $U_{r}$ and $V_{r}^T$ by the square root of $\Sigma_{r}$ to obtain $Y_r$ and $Z_r$: 
		\begin{equation} \label{eq:svd_2_factors}
		X\approx X_{r} = \sum_{i=1}^{r} y_{i}z_{i} = Y_{r} Z_{r}, 
		\end{equation} 
		where $Y_r = U_r \Sigma_r^{1/2}$, $Z_r =  \Sigma_r^{1/2} V_r^T$, $y_i = \sqrt{\sigma_i} u_i$ and $z_i = \sqrt{\sigma_i} v_i^T$ for $1 \leq i \leq r$. 
		Matrices \(Y_{r}\) and \(Z_{r}\) cannot be used directly for NMF initialization since \(Y_{r}\) and \(Z_{r}\) usually contain negative elements (roughly half of them, except for the first factor, by the Perron-Frobenius theorem~\cite{berman1994nonnegative}).

		Given a vector $x$, let us denote $x^{(\geq0)} = \max(0,x)$ its nonnegative part and 
		$x^{(\leq 0)} = \max(0,-x)$  its nonpositive part so that $x = x^{(\geq0)} - x^{(\leq 0)}$. 
		Using this notation, \eqref{eq:svd_2_factors} can be rewritten as: 
		\begin{equation}\label{eq:svd_2_factors_split}
		X\approx X_{r} = \sum_{i=1}^{r} y_{i}z_{i}  = \sum_{i=1}^{r}\Big( y_{i}^{(\geq0)}z_{i}^{(\geq0)}+y_i^{(\leq 0)}z_{i}^{(\leq 0)}\Big)
		\; - \; 
		\sum_{i=1}^{r}\Big(y_{i}^{(\geq0)}z_{i}^{(\leq 0)} + y_{i}^{(\leq 0)}z_{i}^{(\geq0)}\Big). 			
		\end{equation}
		To obtain a feasible initialization for NMF, we have to deal with the second summand which leads to negative elements in the decomposition.  
		Currently, there are mostly two approaches used in practice for this purpose. 
		
		The first approach discards the second summand and selects \(r\) product terms from the first summand on the basis of  some criterion. In particular, the most widely used method, namely nonnegative double SVD (\textsc{NNDSVD})~\cite{boutsidis_svd_2008}, selects \(r\) terms as follows: for each \(i\), it selects  
		$y_{i}^{(\geq0)}z_{i}^{(\geq0)}$ if $||y_{i}^{(\geq0)}z_{i}^{(\geq0)}||_F >  ||y_i^{(\leq 0)}z_{i}^{(\leq 0)}||_F$, 
		otherwise it selects $y_i^{(\leq 0)}z_{i}^{(\leq 0)}$.  
		This is equivalent to projecting $Y_r$ and $Z_r$ onto the nonnegative orthant but taking advantage of the sign ambiguity of the SVD~\cite{bro2008resolving}. 
		The second approach takes the absolute value of the second term, which is equivalent to using  
		$W = |Y_r|$ and $H = |Z_r|$ as an initialization for NMF~\cite{qiao_new_2015}. This method is referred to as SVD-NMF.    
	
	Let us denote $X^{\geq 0}_{r}$ the solution obtained by one of the two approaches mentioned above. In both cases, we will have 
	\[
	X^{\geq 0}_{r+1} \; \geq \; X^{\geq 0}_{r} \quad \text{ for all } r \geq 1, 
	\] 
	 since each rank-one factor selected from the SVD is nonnegative. Hence, for $r$ sufficiently large, the error $||X-X^{\geq 0}_{r}||_F$ will increase as $r$ increases since the negative terms are not taken into account; 
	see Figure~\ref{fig:residual_curves} for examples on real data sets.
		Like the unconstrained rank-$r$ approximation $X_r$ of \(X\), 
		it would make sense that the approximation quality of \(X^{\geq 0}_{r}\) increases as $r$ increases. Another drawback of these approaches is that they 
		either throw away half of the rank-one factors of the first summand and all of the rank-one factors in the second summand (as in NNDSVD) or sum them together so that the sign information is lost (as in SVD-NMF): a lot of information is wasted. 
		
		In order to avoid these two important drawbacks, we propose a new method where 
		\begin{itemize}
		    \item[(i)] We keep all the terms from the first summand in~\eqref{eq:svd_2_factors_split}. 
		Hence, we will only need a truncated SVD of rank $\lceil \frac{r}{2} + 1 \rceil$. 
		In fact, assuming the matrices $XX^T$ and $X^TX$ are irreducible\footnote{A symmetric matrix is irreducible if and only if its associated graph is connected.} (which is the case for all the matrices we have tested in practice), the first rank-one factor $y_1z_1$ of the SVD is positive, by the Perron-Frobenius theorem~\cite{berman1994nonnegative}. 
		This implies that $y_{i}^{(\geq0)} z_{i}^{(\geq 0)} \neq 0$ 
		and $y_{i}^{(\leq 0)} z_{i}^{(\leq 0)} \neq 0$ 
		for all $i \geq 2$ because the singular triplets are orthogonal to one another~\cite{golub_matrix_2013}, that is, $y_i^T y_1 = z_i z_1^T  = 0$ for all $i \geq 2$, which implies that $y_i$ and $z_i$ contain at least one positive and one negative entry. 
		    
		    \item[(ii)] Although we also discard the second summand as in NNDSVD, we will use this information to improve the terms in the first summand. This can be done computationally very efficiently using the low-rank structure of the second summand; see the details below. 
		\end{itemize}

	Our initialization is described in Algorithm~\ref{algo:nnsvdlrc}. 
	It works as follows: Let $p = \lceil r/2 + 1 \rceil$. Then, 
	\begin{enumerate}
	
	\item Compute the rank-$p$ truncated SVD of $X$, with $X_p = \sum_{i=1}^p y_i z_i$; see~\eqref{eq:svd_2_factors}.  
	
	    \item The first rank-one factor of the SVD is used to initialize $W(:,1)$ and $H(1,:)$, that is, 
	    \[
	    W(:,1) = |y_1| 
	    \quad \text{ and } \quad
	    H(1,:) = |z_1|. 
	    \]
	    Note that the absolute value is used because the SVD has a sign ambiguity (hence could generate $y_1$ and $z_1$ with negative entries). In any case, $|y_1| |z_1|$ is an optimal rank-one approximation since $X$ is nonnegative~\cite{berman1994nonnegative}. 
	    
	    \item The other $r-1$ rank-one factors are given by the next $\lceil r/2 \rceil$ factors of the truncated SVD as follows: 
	    \[ 
	    W(:,i) = y_i^{(\geq 0)}, \;  
	    W(:,i+1) = y_i^{(\leq 0)}, \; 
	    H(i,:) = z_i^{(\geq 0)} \, \text{ and } \, 
	    H(i+1,:) = z_i^{(\leq 0)}, 
	    \]
	    where $i = 2, 4, \dots$, in order to obtain a nonnegative NMF initialization $(W,H)$ with $r$ factors. 
	    Note that, by this construction, the average sparsity of these factors is at least 50\%. (In practice, SVD factors usually do not contain zero entries hence average sparsity is exactly 50\%, ignoring the first rank-one factor.) 
	    
	    \item In order to improve the current solution $(W,H)$ built using the first $p$ singular triplets, 
	    we propose to update them using the low-rank approximation $X_p$ by performing a few iteration of an NMF algorithm on the problem 
	    \[
	    \min_{W \geq 0, H \geq 0} ||X_p - WH||_F^2, \quad \text{ where } X_p = Y_p Z_p. 
	    \]
	    The reason for this choice is that, for most NMF algorithms,  performing such iterations is significantly cheaper than performing a standard NMF iteration on the input matrix $X$. 
	    In fact, the most expensive steps of most NMF algorithms is to compute $X H^T$, $W^T X$, $HH^T$ and $W^TW$ which relates to computing the gradient of the objective function; see, e.g., \cite{gillis_accelerated_2012}. 
	    When $X=X_p$ has a low-rank representation $X_p = Y_p Z_p$, the cost of one NMF iteration reduces from $O(mnr)$ operations to $O((m+n)r^2)$ operations.  
	In this paper, we use the state-of-the-art NMF algorithm referred to as accelerated hierarchical alternating least squares (A-HALS)~\cite{gillis_accelerated_2012} to perform this step. A proper implementation requires $O((m+n)r^2)$ operations per iteration instead of $O(mnr)$ if we would apply A-HALS on the input matrix $X$, as explained above. 
	We run A-HALS as long as the relative error decreases the initial error by a proportion of $\delta$. 
	We used $\delta = 5\%$ which leads in all tested cases to less than 10 iterations, which are negligible compared to computing the truncated SVD that requires $\Omega(pmn)$ operations, and to the subsequent NMF iterations, that require $O(mnr)$ operations. 
	    
	    	The idea of using a low-rank approximation of $X$ to speep up NMF computations was proposed in~\cite{zhou_fast_2012}, but not in combination with A-HALS nor as an initialization procedure. 
	\end{enumerate}

		For these reasons, we will refer to our method as nonnegative SVD with low-rank correction (\textsc{NNSVD-LRC}) as it consist of 
		(i)~a selection of nonnegative factors   from the \textsc{SVD} followed by 
		(ii)~NMF iterations that uses the low-rank approximation $X_p$ of $X$, for a negligible additional computational cost of $O((m+n)r^2)$ operations.

\algsetup{indent=2em}
\begin{algorithm}[ht!]
\caption{Nonnegative Singular Value Decomposition with Low-Rank Correction (\textsc{NNSVD-LRC})}
			\label{algo:nnsvdlrc} 
\begin{algorithmic}[1]
\REQUIRE An $m$-by-$n$ nonnegative matrix $X$ and a positive integer $r$. 

    \ENSURE Nonnegative factors $W \in \mathbb{R}^{m \times r}$ and $H \in \mathbb{R}^{r \times n}$ such that $X \approx WH$

	 \STATE		p = $\lceil r/2 + 1 \rceil$;
    
    \STATE	$[U, \Sigma, V]$ = truncated-SVD($X$, $p$);

	 \STATE		$Y_p = U \Sigma^{1/2}$; $Z_p= \Sigma^{1/2} V^{T}$;
			
	  \STATE	\emph{\% Populating \(W\) and \(H\) using $Y_p$ and $Z_p$} 
	 \STATE		$W(:,1) = |Y_p(:, 1)|$;	$H(1,:) = |Z_p(1,:)|$; 
	 
		\STATE $i = 2$; $j = 2$; 
		
		 \WHILE{$i \leq r$}
		 
	\IF{$i$ is even}
       \STATE $W(:,i) = \max(Y_p(:,j),0)$; $H(i,:) = \max(Z_p(j,:),0)$; 
    \ELSE
        \STATE $j = j+1$; 
        \STATE $W(:,i) = \max(-Y_p(:,j),0)$; $H(i,:) = \max(-Z_p(j,:),0)$; 
    \ENDIF

     \STATE $i = i+1$; 
		 
    \ENDWHILE 

		\STATE $e_0 = ||X_p - WH||_F$; $k = 0$; 
			\STATE	\emph{\% Improve \(W\) and \(H\) by applying A-HALS on the low-rank matrix $X_p=Y_pZ_p$}
			
			\WHILE{ $k = 0$ or $e_k - e_{k-1} \geq \delta e_0$ } 
			
			\STATE Perform one iteration of A-HALS on $X_p = Y_pZ_p$ starting from $(W,H)$ to obtain an improved solution $(W,H)$. 
			
			\STATE $e_{k+1} = ||X_p - WH||_F$;
			
			\STATE $k = k+1$; 
			
			\ENDWHILE 

\end{algorithmic}
\end{algorithm}

	\begin{remark}[Computation of the error]
	In Algorithm~\ref{algo:nnsvdlrc}, the error $||X_p-WH||_F$ has to be computed: this can be done in $O((m+n)r^2)$ operations observing that 
	\begin{align*}
	||X_p-WH||_F^2 & = \langle X_p, X_p \rangle - 2 \langle X_p, WH \rangle + \langle WH, WH \rangle \\ 
	               & = \langle Y_p Z_p, Y_p  Z_p \rangle - 2 \langle Y_p Z_p, WH \rangle + \langle W^TW, HH^T \rangle \\ 
	               & = \langle Y_p^T Y_p, Z_p Z_p^T \rangle - 2 \langle (W^T Y_p) Z_p , H  \rangle + \langle W^TW, HH^T \rangle, 
	\end{align*}
	where $\langle A, B \rangle = \sum_{i,j} A_{i,j} B_{i,j}$ is the inner product associated with the Frobenius norm.  
	\end{remark}

	\section{Numerical Experiments} \label{numexp}
	
	In this section, we compare NNSVD-LRC with NNDSVD and SVD-NMF. All tests are preformed using All tests are preformed using Matlab R2017b (Student License) on a laptop Intel CORE i5-2540M CPU @2.60GHz 4GB RAM. 
	The code is available from \url{https://sites.google.com/site/nicolasgillis/code}.  
		Due to the space limit, we restrict ourselves to three dense and three sparse widely used data sets; see Tables~\ref{tb:biometric_datasets} and~\ref{tb:text_mining_datasets}. 
		We also restrict ourselves to using the multiplicative update algorithm, one of the most widely used one. 
		(On the Matlab code provided online, we provide experiments for two other data sets, namely the CBCL facial images, and the classic document data set, in combination with A-HALS.) 
		\begin{table}[!h]
			\caption{Biometric data sets}
			\label{tb:biometric_datasets}
			\begin{center}
			\begin{threeparttable}
				\begin{tabular}{|c | c | c | c |} 
					\hline
					Data set & Image size (\(h\times w\)) & \(m = h\times w\) & \(n\) \\ 
					\hline\hline
					AT\&T Faces$^{a}$~\cite{qiao_new_2015} & \(112\times 92\) & 10304 & 400 \\ 
					\hline
					IITD Iris$^{b}$~\cite{kumar_comparison_2010} & \(240\times 320\) & 76800 & 200 \\
					\hline
					TD Fingerprints$^{c}$~\cite{si_detection_2015} & \(750 \times 800\) & 600000 & 100 \\
					\hline
				\end{tabular}
				\begin{flushleft}
				\footnotesize  
				
				$^{a}$\url{http://www.cl.cam.ac.uk/research/dtg/attarchive/facedatabase.html}
				
				$^{b}$\url{http://www4.comp.polyu.edu.hk/~csajaykr/IITD/Database_Iris.htm}
				
				$^{c}$\url{http://ivg.au.tsinghua.edu.cn/dataset/TDFD.php}
				\end{flushleft}
			\end{threeparttable}
			\end{center}
		\end{table}
		\begin{table}[!h]
			\caption{Document data sets from~\cite{ZG05}}
			\label{tb:text_mining_datasets}
			\begin{center}
				\begin{tabular}{|c | c | c | c | c |} 
					\hline
					Dataset Name & \#nonzeros & Sparsity & $m$ & $n$ \\ 
					\hline\hline
					Sports & 1091723 & 99.14 & 8580 & 14870 \\ 
					\hline
					Reviews & 758635 & 98.99 & 4069 & 18483 \\
					\hline
					Hitech & 331373 & 98.57 & 2301 & 10080 \\
					\hline
				\end{tabular}
			\end{center}
		\end{table}

	Throughout this section, we will use the following two quantities: 
	\begin{enumerate}
		\item the relative error which measures the quality of an NMF solution: 
		\begin{equation*}
		\text{relative error}(W,H) \; = \;  \frac{\|X - WH\|_{F}}{\|X\|_{F}}, 
		\end{equation*}
		\item the sparsity which measures the proportion of zero entries in a matrix:
		\begin{equation*}
		\text{sparsity}(W) = \frac{\text{\# of zeros in $W$}}{\text{\# of total elements in $W$}}.  
		\end{equation*}
	\end{enumerate}

\begin{figure}[h!] 
		\centering
			\begin{tabular}{cc} 
		 \includegraphics[width=0.48\textwidth]{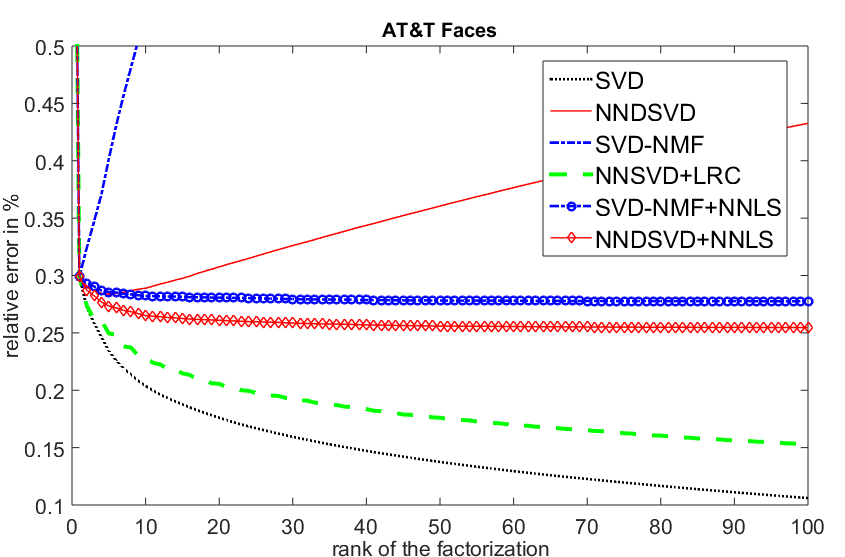}  &  
		\includegraphics[width=0.48\textwidth]{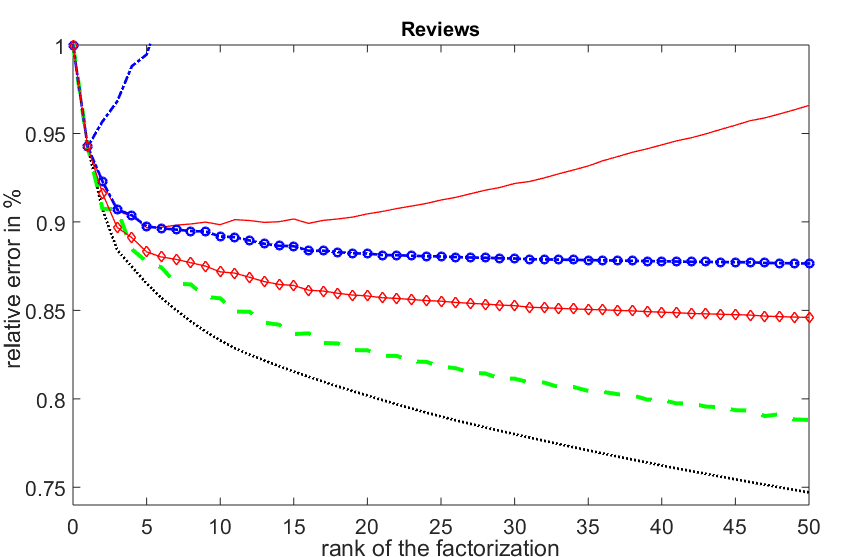}  \\
				\includegraphics[width=0.48\textwidth]{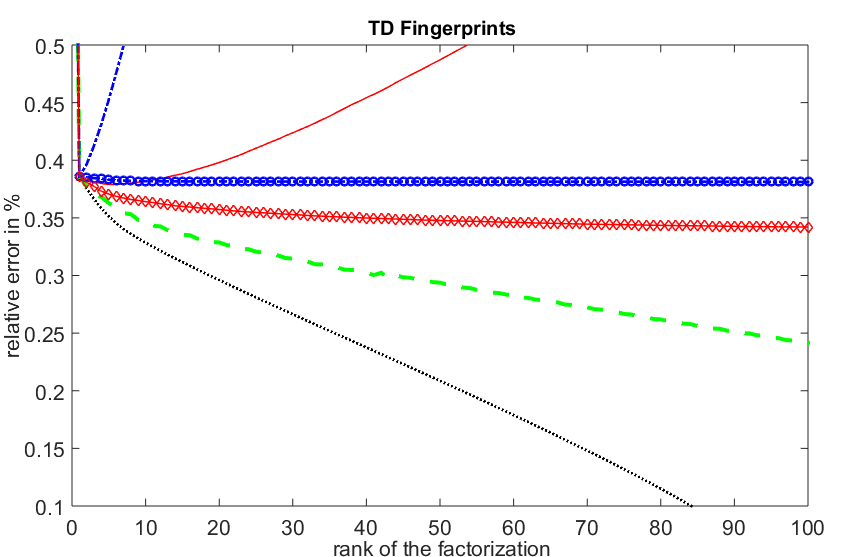} & 
					\includegraphics[width=0.48\textwidth]{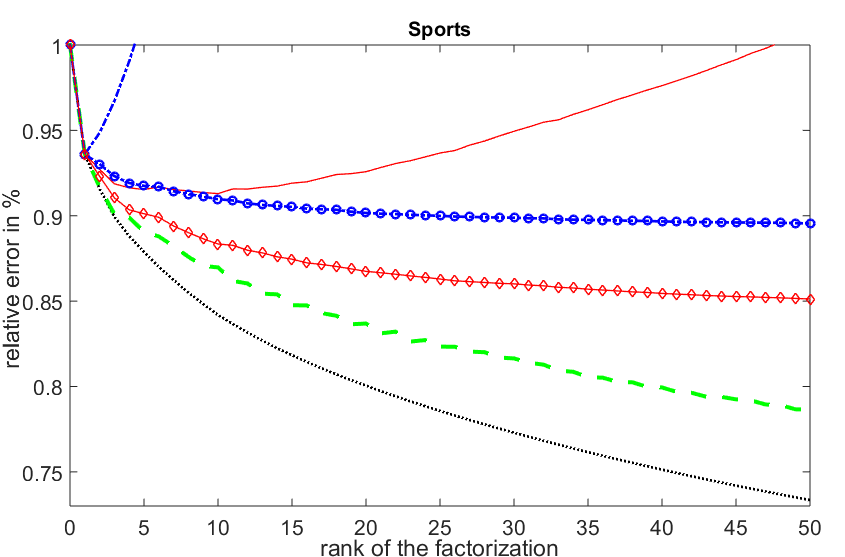} \\  
		  \includegraphics[width=0.48\textwidth]{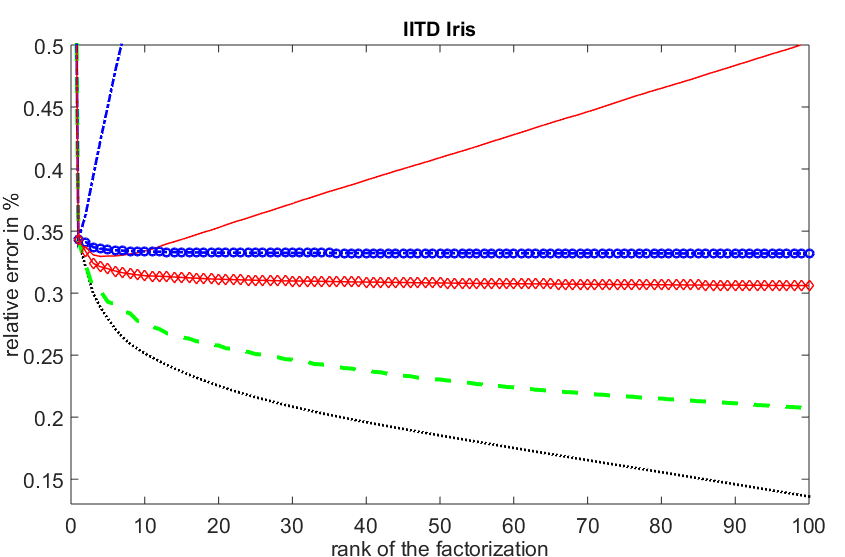} & 
		\includegraphics[width=0.48\textwidth]{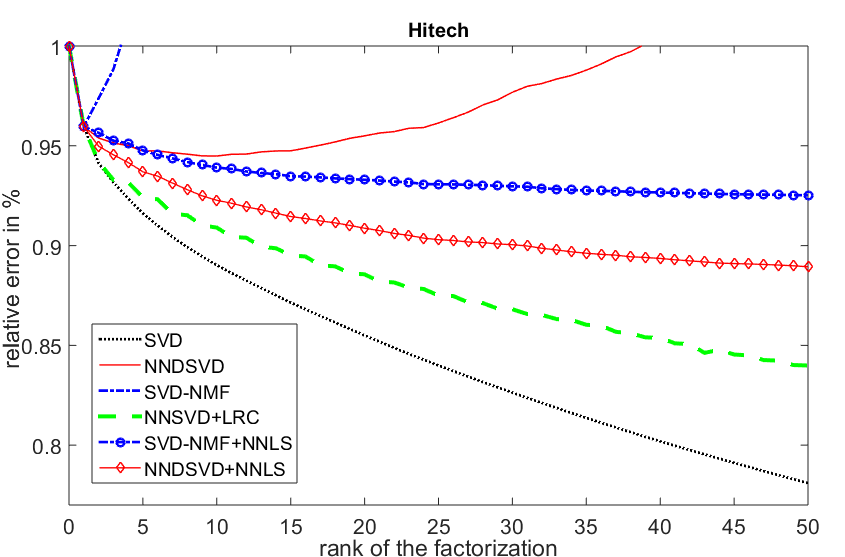} \\ 
		\end{tabular}
			\caption{Relative error of the SVD-based NMF initializations for different values of the rank $r$. \label{fig:residual_curves}} 
		\end{figure}

	\paragraph{Initial error} 		Figure~\ref{fig:residual_curves} displays the relative errors in percent for different values of $r$ for each data set.  
		This illustrates the fact that the error of \textsc{NNDSVD} and \textsc{SVD-NMF} increases as $r$ increases (as soon as $r$ is sufficiently large); 
		see the discussion in Section~\ref{proposedsol}.  
		In contrast, the error of \textsc{NNSVD-LRC} decreases as $r$ increases. 
		Note that the relative error of \textsc{SVD-NMF} grows much faster than \textsc{NNDSVD}.    

		One may argue that the above comparison is not totally fair as \textsc{SVD-NMF} and \textsc{NNDSVD} did not update the factors $W$ and $H$ as opposed to \textsc{NNSVD-LRC}. 
		Therefore, Figure~\ref{fig:residual_curves} also displays the relative error of these initializations after the matrix \(H\) is updated with the solution of the nonnegative least squares (NNLS) problem $\min_{H \geq 0} ||X-WH||_F$ for $W$ fixed. 
			This allows to compare the quality of the basis matrix generated by the different initializations. 
		We observe that \textsc{NNSVD-LRC} still outperforms \textsc{SVD-NMF} and \textsc{NNDSVD} after this update. 
		Table~\ref{tb:superior} displays the relative error in percent of the three SVD-based initializations for different values of the factorization rank $r$, after the NNLS update, and also after one iteration of the HALS algorithm. Although the error of \textsc{SVD-NMF} and \textsc{NNDSVD} decreases significantly compared to the initial error (cf.\@ Figure~\ref{fig:residual_curves}), it is still much higher than \textsc{NNSVD-LRC}.  
		\begin{table}[!h]
			\centering
			\caption{Comparison of the relative error (in percent) of the SVD-based NMF initializations when they are aided by one iteration of \textsc{HALS} and the \textsc{NNLS} update of $H$. The lowest error in highlighted in bold.}
			\label{tb:superior}
			\begin{tabular}{|c||c|c|c||c|c|c||c|c|c|}
				\hline
				\multirow{2}{*}{ } & \multicolumn{3}{c||}{AT\&T} & \multicolumn{3}{c||}{IITD} & \multicolumn{3}{c|}{TDF} \\ \cline{2-10} 
			 	& $r$=60  & $r$=80 & $r$=100 & $r$=30  & $r$=40 & $r$=50 & $r$=15  & $r$=20 & $r$=25 \\ \hline \hline
				NNDSVD+HALS    & 22.10   & 21.71  & 21.35   & 27.69   & 27.28  & 26.94  & 35.03   & 34.58  & 34.22 \\ \hline
				NNDSVD+NNLS     & 25.55   & 25.49  & 25.46   & 30.98   & 30.89  & 30.82  & 35.99   & 35.73  & 35.47 \\ \hline
				SVD-NMF+HALS   & 22.14   & 21.37  & 20.76   & 28.61     & 27.85  & 27.28  & 36.47   & 36.17  & 35.81 \\ \hline
				SVD-NMF+NNLS    & 27.80   & 27.77  & 27.76   & 33.24   & 33.23  & 33.22  & 38.17   & 38.16  & 38.16 \\ \hline
				NNSVD-LRC     & \textbf{17.00}   & \textbf{16.05}  & \textbf{15.29}   & \textbf{24.63}   & \textbf{23.76}  & \textbf{23.04}  & \textbf{33.54}   & \textbf{32.87}  & \textbf{32.06} \\ \hline
				\multirow{2}{*}{ } & \multicolumn{3}{c||}{Sports} & \multicolumn{3}{c||}{Reviews} & \multicolumn{3}{c|}{Hitech} \\ \cline{2-10} 
				 & $r$=15  & $r$=20  & $r$=25  & $r$=15   & $r$=20  & $r$=25  & $r$=15   & $r$=20  & $r$=25  \\ \hline \hline
				NNDSVD+HALS    & 85.69   & 84.65   & 83.90   & 84.19    & 83.35   & 82.78   & 89.93    & 89.09   & 88.29  \\ \hline
				NNDSVD+NNLS     & 87.46   & 86.74   & 86.29   & 86.41    & 85.82   & 85.52   & 91.46    & 90.87   & 90.29  \\ \hline
				SVD-NMF+HALS   & 87.04   & 86.12   & 85.46   & 84.83    & 84.10   & 83.64   & 90.72    & 90.02   & 89.45  \\ \hline
				SVD-NMF+NNLS    & 90.52   & 90.16   & 89.99   & 88.63    & 88.21   & 88.05   & 93.48    & 93.29   & 93.07  \\ \hline
				NNSVD-LRC     & \textbf{84.76}   & \textbf{83.69}   & \textbf{82.33}   & \textbf{83.66}    & \textbf{82.74}   & \textbf{81.82}   & \textbf{89.51}    & \textbf{88.56}   & \textbf{87.50}  \\ \hline
			\end{tabular}
		\end{table}

		\paragraph{Sparsity} For the sparsity of the initializations, \textsc{SVD-NMF} generates dense initial factors, with sparsity 0\% in all cases (because SVD generated dense factors and SVD-NMF take their absolute values as initial estimates for $W$ and $H$).  
		 \textsc{NNDSVD} generates factors with average sparsity 49\%, 
		with the sparsity of every initialization $(W,H)$ being between 45\% and 53\% for all data sets. 
		\textsc{NNSVD-LRC} generates factors with average sparsity 45\% (resp.\@ 58\%), 
		with the sparsity of every initialization $(W,H)$ being between 23\% (resp.\@ 51\%) and 59\% (resp.\@ 66\%) for dense (resp.\@ sparse) data sets. 
		This confirms our discussion in Section~\ref{proposedsol} where the initialization provided by \textsc{NNSVD-LRC} has average sparsity around 50\%, similarly as 
		\textsc{NNDSVD}. (Note that this is not exactly 50\% because of the low-rank correction step performed by \textsc{NNSVD-LRC}.) 

	\paragraph{Computational time} 
	
	Table~\ref{tb:cpu_time} reports the computational time for the different initializations on the different data sets, averaged over 100 runs. 
	As expected, \textsc{NNDSVD} and \textsc{SVD-NMF} have roughly the same computational cost, the main cost being the computation of the rank-$r$ truncated SVD, while 
	\textsc{NNSVD-LRC} is faster as the main computational cost is the computation of the rank-$p$ truncated SVD, with $p = \lceil r/2 + 1 \rceil$, with an additional cost of running A-HALS on the rank-$p$ approximation of $X$. 
	
	\begin{table}[!h]
			\centering
			\caption{CPU time (in s.) taken by different \textsc{NMF} initializations for the  different data sets. Bold indicates the algorithm that took less CPU time.}  
			\label{tb:cpu_time}
			\begin{tabular}{|l||c|c|c||c|c|c||c|c|c|}
				\hline
				\multirow{2}{*}{}  		& \multicolumn{3}{c||}{AT\&T} & \multicolumn{3}{c||}{IITD} & \multicolumn{3}{c|}{TDF} \\ \cline{2-10} 
			 	& $r$=60     & $r$=80     & $r$=100     & $r$=30     & $r$=40      & $r$=50      & $r$=15     & $r$=20     & $r$=25 \\ \hline \hline
				NNDSVD    & 4.10    & 5.68   & 7.72    & 28.59    & 65.97   & 84.89  & 18.65   & 14.91   & 15.66   \\ \hline
				SVD-NMF       & 4.09    & 5.70   & 7.72    & 28.46   & 65.53  & 84.48   & 18.54   & 14.83   & \textbf{15.64}   \\ \hline
				NNSVD-LRC          & \textbf{2.77}    & \textbf{3.89}   & \textbf{5.19}    & \textbf{26.99}    & \textbf{48.11}   & \textbf{60.09}   & \textbf{14.56}   & \textbf{14.71}   & 15.68   \\ \hline
				\hline
				\multirow{2}{*}{ } & \multicolumn{3}{c||}{Sports} & \multicolumn{3}{c||}{Reviews} & \multicolumn{3}{c|}{Hitech} \\ \cline{2-10} 
			 & $r$=15 & $r$=20 & $r$=25 & $r$=15 & $r$=20 & $r$=25 & $r$=15 & $r$=20 & $r$=25 \\ \hline \hline
				NNDSVD    & 4.30   & 4.63  & 5.96  & 3.13   & 3.39    & 4.40  & 1.89   & 2.17   & 2.83   \\ \hline
				SVD-NMF       & 4.31   & 4.62   & 5.95   & 3.10   & 3.39    & 4.40   & 1.89   & 2.17   & 2.82   \\ \hline
				NNSVD-LRC          & \textbf{3.19}    & \textbf{3.83}    & \textbf{5.05}  & \textbf{2.52}    & \textbf{2.97}     & \textbf{4.03}    & \textbf{1.54}    & \textbf{1.93}    & \textbf{2.27}   \\ \hline
			\end{tabular}
		\end{table}


	\paragraph{Convergence of NMF algorithms} 
	
	We now compare the three NMF initializations used in combination with one of the most widely used NMF algorithm, namely, the multiplicative updates (MU)~\cite{lee_learning_1999, lee_algorithms_2001}. 
	Table~\ref{tb:convergenceMU} displays the relative error in percent after 1, 10 and 100 iterations of MU. 
		\begin{table}[!h]
			\centering
			\caption{Relative error in percent of MU after 1, 10 and 100 iterations when seeded by different SVD-based NMF initializations on the dense and sparse data sets. 
			The lowest error is highlighted in bold.} 
			\label{tb:convergenceMU}
			\begin{tabular}{|c|c||c|c|c||c|c|c||c|c|c|}
				\hline
				\multicolumn{2}{|l||}{\multirow{2}{*}{  } } & \multicolumn{3}{c||}{AT\&T} & \multicolumn{3}{c||}{IITD} & \multicolumn{3}{c|}{TDF} \\ \cline{3-11} 
				\multicolumn{2}{|r||}{ }                
				& $r$=60  & $r$=80 & $r$=100 & $r$=30  & $r$=40 & $r$=50 & $r$=15 & $r$=20 & $r$=25 \\ \hline \hline
				NNDSVD        & \multirow{3}{*}{1} & 24.58  & 24.51 & 24.47  & 29.73  & 29.56 & 29.44 & 35.88 & 35.55 & 35.26 \\ \cline{1-1} \cline{3-11} 
				SVD-NMF       &                         & 30.03  & 30.02 & 30.02  & 34.40  & 34.41 & 34.40 & 38.97 & 39.06 & 39.12 \\ \cline{1-1} \cline{3-11} 
				NNSVD-LRC &                         & \textbf{16.91}  & \textbf{15.96} & \textbf{15.19}  & \textbf{24.53}  & \textbf{23.66} & \textbf{22.94} & \textbf{33.49} & \textbf{32.81} & \textbf{31.99} \\ \hline \hline
				NNDSVD        & \multirow{3}{*}{10} & 21.71  & 21.52 & 21.40  & 26.99  & 26.68 & 26.43 & 34.43 & 33.81 & 33.27 \\ \cline{1-1} \cline{3-11} 
				SVD-NMF       &                         & 27.18  & 27.15 & 27.14  & 31.66  & 31.60 & 31.55 & 37.11 & 36.96 & 36.85 \\ \cline{1-1} \cline{3-11} 
				NNSVD-LRC &                         & \textbf{16.62}  & \textbf{15.67} & \textbf{14.91}  & \textbf{24.19}  & \textbf{23.33} & \textbf{22.61} & \textbf{33.35} & \textbf{32.62} & \textbf{31.77} \\ \hline \hline
				NNDSVD        & \multirow{3}{*}{100} & 17.83  & 17.09 & 16.52  & 24.40  & 23.69 & 23.13 & 33.37 & 32.42 & 31.55 \\ \cline{1-1} \cline{3-11} 
				SVD-NMF       &                         & 17.06  & 16.40 & 15.92  & \textbf{23.52}  & 22.69 & 22.10 & \textbf{32.38} & \textbf{31.18} & \textbf{30.06} \\ 
				\cline{1-1} \cline{3-11} 
				NNSVD-LRC &                         & \textbf{15.94}  & \textbf{14.96} & \textbf{14.19}  & 23.54  & \textbf{22.64} & \textbf{21.94} & 32.91 & 32.05 & 31.10 \\ \hline  
				\hline
				\multicolumn{2}{|c||}{\multirow{2}{*}{}} & \multicolumn{3}{c||}{Sports} & \multicolumn{3}{c||}{Reviews} & \multicolumn{3}{c|}{Hitech} \\ \cline{3-11} 
				\multicolumn{2}{|r||}{rank $r$}           
				       & $r$=15 & $r$=20 & $r$=25 & $r$=15 & $r$=20 & $r$=30 & $r$=15 & $r$=20 & $r$=25  \\ \hline \hline
				NNDSVD        & \multirow{3}{*}{1} & 87.22  & 86.53  & 86.07  & 85.57   & 84.91  & 84.58  & 91.13  & 90.56  & 89.99  \\ \cline{1-1} \cline{3-11} 
				SVD-NMF       &                         & 90.90  & 90.56  & 90.46  & 88.89   & 88.54  & 88.38  & 93.60  & 93.48  & 93.31  \\ \cline{1-1} \cline{3-11} 
				NNSVD-LRC &                         & \textbf{84.42}  & \textbf{83.33}  & \textbf{82.00}  & \textbf{83.32}   & \textbf{82.51}  & \textbf{81.55}  & \textbf{89.29}  & \textbf{88.18}  & \textbf{86.99}  \\ \hline \hline
				NNDSVD        & \multirow{3}{*}{10} & 84.17  & 82.70  & 81.40  & 82.78   & 81.66  & 81.01  & 88.48  & 87.34  & 86.24  \\ \cline{1-1} \cline{3-11} 
				SVD-NMF       &                         & 84.02  & 82.68  & 81.34  & 83.01   & 81.90  & 81.03  & \textbf{88.35}  & 86.89  & \textbf{85.59}  \\ \cline{1-1} \cline{3-11} 
				NNSVD-LRC &                         & \textbf{83.67}  & \textbf{81.80}  & \textbf{80.70}  & \textbf{82.42}   & \textbf{81.54}  & \textbf{80.45}  & 88.40  & \textbf{86.83}  & 85.63  \\ \hline \hline
				NNDSVD        & \multirow{3}{*}{100} & \textbf{82.93}  & 81.38  & 79.89  & 82.30  & 80.96  & 79.92  & 87.88  & 86.49  & 85.10  \\ \cline{1-1} \cline{3-11} 
				SVD-NMF       &                         & 82.96  & 81.09  & \textbf{79.55}  & 82.20   & \textbf{80.76}  & \textbf{79.63}  & \textbf{87.70}  & \textbf{86.14}  & \textbf{84.79}  \\ \cline{1-1} \cline{3-11} 
				NNSVD-LRC &                         & 83.00  & \textbf{81.05}  & 79.85  & \textbf{82.09}   & 81.05  & 79.79  & 87.94  & 86.34  & 85.02  \\ 		\hline
			\end{tabular}
		\end{table}

	We observe the following: 
	\begin{itemize}
	\item \textsc{NNDSVD} and \textsc{SVD-NMF} with 1 or 10 iterations of MU are not enough to get back at \textsc{NNSVD-LRC}, except for the Hitech data set where \textsc{SVD-NMF} achieves a slightly lower error (0.05\% for $r=15$ and 0.04\% for $r=25$). 
	This is explained by the fact that the inital error of \textsc{NNSVD-LRC} is much lower, as shown in Figure~\ref{fig:residual_curves} and Table~\ref{tb:superior}. 
	
	    \item After 100 iterations of the MU, \textsc{NNDSVD} and \textsc{SVD-NMF} sometimes are able to get back at \textsc{NNSVD-LRC}: there is no clear winner 
	    (although on these 6 data sets, \textsc{NNDSVD} seems to perform worse).
	    The MU have converged (close) to different stationary points and there is no guarantee in general that \textsc{NNSVD-LRC} will lead to better local solutions. 
	\end{itemize}
	
	
In summary, \textsc{NNSVD-LRC} is able to obtain a better (and sparse) initial solution faster than \textsc{NNSDVD} and \textsc{SVD-NMF}.
It should therefore always be preferred if one wants to quickly obtain a good solution. 
However, due to the complexity of NMF~\cite{vavasis_complexity_2010}, 
if one wants to obtain a possibly better solution, 
it is recommended to use multiple initializations and keep the best solution obtained; see, e.g., \cite{cichocki2009nonnegative} for a discussion.  
	
	\section{Conclusion} \label{concl} 
	
	In this paper, we presented a novel \textsc{SVD}-based \textsc{NMF} initialization. Our motivation was to address the shortcomings of previously proposed \textsc{SVD}-based \textsc{NMF} initializations. Our newly proposed method, referred to as nonnegative singular value decomposition with low-rank correction (NNSVD-LRC), has the following advantages 
	\begin{enumerate}
	    \item the initial error decreases as the factorization $r$ increases, 
	    
	    \item the average sparsity of the initial factors $(W,H)$ is close to 50\%, 
	    
	    \item it is computationally cheaper as it only requires the computation of a truncated SVD of rank $p=\lceil r/2+1\rceil$, instead of $r$, and 
	    
	    \item it takes advantage of the discarded factors using highly efficient NMF iterations based on the low-rank approximation computed by the SVD.

	\end{enumerate}
In summary, NNSVD-LRC provides better initial NMF factors (both in terms of error and sparsity) at a lower computational cost. This was confirmed on both dense and sparse real data sets. 
This allows NMF algorithms to converge faster to a stationary point, 
although there is no guarantee that this stationary point will have lower error than other initializations, 
as NMF is a difficult non-convex optimization problem~\cite{vavasis_complexity_2010}.


	\section*{Acknowledgement} 
		The financial support of HEC Pakistan is highly acknowledged for granting PhD scholarship to the first author. 
		NG acknowledges the support of the European Research Council (ERC starting grant n$^\text{o}$ 679515). 

\small 

\bibliographystyle{spmpsci} 
\bibliography{NNSVD-LRC_paper}

\begin{thebibliography}{10}
\providecommand{\url}[1]{{#1}}
\providecommand{\urlprefix}{URL }
\expandafter\ifx\csname urlstyle\endcsname\relax
  \providecommand{\doi}[1]{DOI~\discretionary{}{}{}#1}\else
  \providecommand{\doi}{DOI~\discretionary{}{}{}\begingroup
  \urlstyle{rm}\Url}\fi

\bibitem{berman1994nonnegative}
Berman, A., Plemmons, R.J.: Nonnegative matrices in the mathematical sciences,
  vol.~9.
\newblock Siam (1994)

\bibitem{boutsidis_svd_2008}
Boutsidis, C., Gallopoulos, E.: {SVD} based initialization: {A} head start for
  nonnegative matrix factorization.
\newblock Pattern Recognition \textbf{41}(4), 1350--1362 (2008).
\newblock \doi{10.1016/j.patcog.2007.09.010}.
\newblock
  \urlprefix\url{http://linkinghub.elsevier.com/retrieve/pii/S0031320307004359}

\bibitem{bro2008resolving}
Bro, R., Acar, E., Kolda, T.G.: Resolving the sign ambiguity in the singular
  value decomposition.
\newblock Journal of Chemometrics \textbf{22}(2), 135--140 (2008)

\bibitem{casalino_subtractive_2014}
Casalino, G., Del~Buono, N., Mencar, C.: Subtractive clustering for seeding
  non-negative matrix factorizations.
\newblock Information Sciences \textbf{257}, 369--387 (2014).
\newblock \doi{10.1016/j.ins.2013.05.038}.
\newblock
  \urlprefix\url{http://linkinghub.elsevier.com/retrieve/pii/S0020025513004349}

\bibitem{cichocki2009nonnegative}
Cichocki, A., Zdunek, R., Phan, A.H., Amari, S.i.: Nonnegative matrix and
  tensor factorizations: applications to exploratory multi-way data analysis
  and blind source separation.
\newblock John Wiley \& Sons (2009)

\bibitem{du_hybrid_2017}
Du, R., Drake, B., Park, H.: Hybrid clustering based on content and connection
  structure using joint nonnegative matrix factorization.
\newblock Journal of Global Optimization  (2017).
\newblock \doi{10.1007/s10898-017-0578-x}.
\newblock \urlprefix\url{http://link.springer.com/10.1007/s10898-017-0578-x}

\bibitem{elad_role_2010}
Elad, M., Figueiredo, M.A.T., {Yi Ma}: On the {Role} of {Sparse} and
  {Redundant} {Representations} in {Image} {Processing}.
\newblock Proceedings of the IEEE \textbf{98}(6), 972--982 (2010).
\newblock \doi{10.1109/JPROC.2009.2037655}.
\newblock \urlprefix\url{http://ieeexplore.ieee.org/document/5420029/}

\bibitem{ensari_character_2016}
Ensari, T.: Character {Recognition} {Analysis} with {Nonnegative} {Matrix}
  {Factorization}.
\newblock International Journal of Computers \textbf{1}, 219--222 (2016)

\bibitem{gillis2014}
Gillis, N.: The why and how of nonnegative matrix factorization.
\newblock Regularization, Optimization, Kernels, and Support Vector Machines
  \textbf{12}(257) (2014)

\bibitem{gillis_using_2010}
Gillis, N., Glineur, F.: Using underapproximations for sparse nonnegative
  matrix factorization.
\newblock Pattern Recognition \textbf{43}(4), 1676--1687 (2010).
\newblock \doi{10.1016/j.patcog.2009.11.013}.
\newblock
  \urlprefix\url{http://linkinghub.elsevier.com/retrieve/pii/S0031320309004324}

\bibitem{gillis_accelerated_2012}
Gillis, N., Glineur, F.: Accelerated {Multiplicative} {Updates} and
  {Hierarchical} {ALS} {Algorithms} for {Nonnegative} {Matrix} {Factorization}.
\newblock Neural Computation \textbf{24}(4), 1085--1105 (2012).
\newblock \doi{10.1162/NECO_a_00256}.
\newblock
  \urlprefix\url{http://www.mitpressjournals.org/doi/10.1162/NECO_a_00256}

\bibitem{golub_matrix_2013}
Golub, G.H., Van~Loan, C.F.: Matrix computations, fourth edition edn.
\newblock Johns {Hopkins} studies in the mathematical sciences. The Johns
  Hopkins University Press, Baltimore (2013).
\newblock OCLC: ocn824733531

\bibitem{tan_using_2011}
Janecek, A., Tan, Y.: Using {Population} {Based} {Algorithms} for
  {Initializing} {Nonnegative} {Matrix} {Factorization}.
\newblock In: Y.~Tan, Y.~Shi, Y.~Chai, G.~Wang (eds.) Advances in {Swarm}
  {Intelligence}, vol. 6729, pp. 307--316. Springer Berlin Heidelberg, Berlin,
  Heidelberg (2011).
\newblock \doi{10.1007/978-3-642-21524-7_37}.
\newblock \urlprefix\url{http://link.springer.com/10.1007/978-3-642-21524-7_37}

\bibitem{kumar_comparison_2010}
Kumar, A., Passi, A.: Comparison and combination of iris matchers for reliable
  personal authentication.
\newblock Pattern Recognition \textbf{43}(3), 1016--1026 (2010).
\newblock \doi{10.1016/j.patcog.2009.08.016}.
\newblock
  \urlprefix\url{http://linkinghub.elsevier.com/retrieve/pii/S0031320309003343}

\bibitem{langville2006initializations}
Langville, A.N., Meyer, C.D., Albright, R., Cox, J., Duling, D.:
  Initializations for the nonnegative matrix factorization.
\newblock In: Proceedings of the twelfth ACM SIGKDD international conference on
  knowledge discovery and data mining, pp. 23--26. Citeseer (2006)

\bibitem{lee_learning_1999}
Lee, D.D., Seung, H.S.: Learning the parts of objects by non-negative matrix
  factorization.
\newblock Nature \textbf{401}, 788 (1999).
\newblock \urlprefix\url{http://dx.doi.org/10.1038/44565}

\bibitem{lee_algorithms_2001}
Lee, D.D., Seung, H.S.: Algorithms for {Non}-negative {Matrix} {Factorization}.
\newblock In: T.K. Leen, T.G. Dietterich, V.~Tresp (eds.) Advances in {Neural}
  {Information} {Processing} {Systems} 13, pp. 556--562. MIT Press (2001).
\newblock
  \urlprefix\url{http://papers.nips.cc/paper/1861-algorithms-for-non-negative-matrix-factorization.pdf}

\bibitem{luce_using_2016}
Luce, R., Hildebrandt, P., Kuhlmann, U., Liesen, J.: Using {Separable}
  {Nonnegative} {Matrix} {Factorization} {Techniques} for the {Analysis} of
  {Time}-{Resolved} {Raman} {Spectra}.
\newblock Applied Spectroscopy \textbf{70}(9), 1464--1475 (2016).
\newblock \doi{10.1177/0003702816662600}.
\newblock
  \urlprefix\url{http://journals.sagepub.com/doi/10.1177/0003702816662600}

\bibitem{maruyama_detecting_2014}
Maruyama, R., Maeda, K., Moroda, H., Kato, I., Inoue, M., Miyakawa, H.,
  Aonishi, T.: Detecting cells using non-negative matrix factorization on
  calcium imaging data.
\newblock Neural Networks \textbf{55}, 11--19 (2014).
\newblock \doi{10.1016/j.neunet.2014.03.007}.
\newblock
  \urlprefix\url{http://linkinghub.elsevier.com/retrieve/pii/S0893608014000707}

\bibitem{prajapati_brain_2015}
Prajapati, S.J., Jadhav, K.R.: Brain tumor detection by various image
  segmentation techniques with introduction to non negative matrix
  factorization.
\newblock Brain \textbf{4}(3), 600--3 (2015)

\bibitem{qiao_new_2015}
Qiao, H.: New {SVD} based initialization strategy for non-negative matrix
  factorization.
\newblock Pattern Recognition Letters \textbf{63}, 71--77 (2015).
\newblock \doi{10.1016/j.patrec.2015.05.019}.
\newblock
  \urlprefix\url{http://linkinghub.elsevier.com/retrieve/pii/S0167865515001762}

\bibitem{rezaei_efficient_2011}
Rezaei, M., Boostani, R., Rezaei, M.: An {Efficient} {Initialization} {Method}
  for {Nonnegative} {Matrix} {Factorization}.
\newblock Journal of Applied Sciences \textbf{11}(2), 354--359 (2011).
\newblock \doi{10.3923/jas.2011.354.359}.
\newblock
  \urlprefix\url{http://www.scialert.net/abstract/?doi=jas.2011.354.359}

\bibitem{sauwen_successive_2017}
Sauwen, N., Acou, M., Bharath, H.N., Sima, D.M., Veraart, J., Maes, F.,
  Himmelreich, U., Achten, E., Van~Huffel, S.: The successive projection
  algorithm as an initialization method for brain tumor segmentation using
  non-negative matrix factorization.
\newblock PLOS ONE \textbf{12}(8), e0180,268 (2017).
\newblock \doi{10.1371/journal.pone.0180268}.
\newblock \urlprefix\url{http://dx.plos.org/10.1371/journal.pone.0180268}

\bibitem{shiga_sparse_2016}
Shiga, M., Tatsumi, K., Muto, S., Tsuda, K., Yamamoto, Y., Mori, T., Tanji, T.:
  Sparse modeling of {EELS} and {EDX} spectral imaging data by nonnegative
  matrix factorization.
\newblock Ultramicroscopy \textbf{170}, 43--59 (2016).
\newblock \doi{10.1016/j.ultramic.2016.08.006}.
\newblock
  \urlprefix\url{http://linkinghub.elsevier.com/retrieve/pii/S0304399116301267}

\bibitem{si_detection_2015}
Si, X., Feng, J., Zhou, J., Luo, Y.: Detection and {Rectification} of
  {Distorted} {Fingerprints}.
\newblock IEEE Transactions on Pattern Analysis and Machine Intelligence
  \textbf{37}(3), 555--568 (2015).
\newblock \doi{10.1109/TPAMI.2014.2345403}.
\newblock \urlprefix\url{http://ieeexplore.ieee.org/document/7029762/}

\bibitem{sun_graph_2016}
Sun, F., Xu, M., Hu, X., Jiang, X.: Graph regularized and sparse nonnegative
  matrix factorization with hard constraints for data representation.
\newblock Neurocomputing \textbf{173}, 233--244 (2016).
\newblock \doi{10.1016/j.neucom.2015.01.103}.
\newblock
  \urlprefix\url{http://linkinghub.elsevier.com/retrieve/pii/S092523121501276X}

\bibitem{vavasis_complexity_2010}
Vavasis, S.A.: On the {Complexity} of {Nonnegative} {Matrix} {Factorization}.
\newblock SIAM Journal on Optimization \textbf{20}(3), 1364--1377 (2010).
\newblock \doi{10.1137/070709967}.
\newblock \urlprefix\url{http://epubs.siam.org/doi/10.1137/070709967}

\bibitem{wang_effective_2012}
Wang, X., Xie, X., Lu, L.: An {Effective} {Initialization} for {Orthogonal}
  {Nonnegative} {Matrix} {Factorization}.
\newblock Journal of Computational Mathematics \textbf{30}(1), 34--46 (2012).
\newblock \doi{10.4208/jcm.1110-m11si10}.
\newblock
  \urlprefix\url{http://www.global-sci.org/jcm/volumes/v30n1/pdf/301-34.pdf}

\bibitem{wild2004improving}
Wild, S., Curry, J., Dougherty, A.: Improving non-negative matrix
  factorizations through structured initialization.
\newblock Pattern recognition \textbf{37}(11), 2217--2232 (2004)

\bibitem{ye_multitask_2015}
Ye, M., Qian, Y., Zhou, J.: Multitask {Sparse} {Nonnegative} {Matrix}
  {Factorization} for {Joint} {Spectral}–{Spatial} {Hyperspectral} {Imagery}
  {Denoising}.
\newblock IEEE Transactions on Geoscience and Remote Sensing \textbf{53}(5),
  2621--2639 (2015).
\newblock \doi{10.1109/TGRS.2014.2363101}.
\newblock \urlprefix\url{http://ieeexplore.ieee.org/document/6939673/}

\bibitem{yoshii_students_2016}
Yoshii, K., Itoyama, K., Goto, M.: Student's {T} nonnegative matrix
  factorization and positive semidefinite tensor factorization for
  single-channel audio source separation.
\newblock pp. 51--55. IEEE (2016).
\newblock \doi{10.1109/ICASSP.2016.7471635}.
\newblock \urlprefix\url{http://ieeexplore.ieee.org/document/7471635/}

\bibitem{ZG05}
Zhong, S., Ghosh, J.: Generative model-based document clustering: a comparative
  study.
\newblock Knowledge and Information Systems \textbf{8 (3)}, 374--384 (2005)

\bibitem{zhou_fast_2012}
Zhou, G., Cichocki, A., Xie, S.: Fast {Nonnegative} {Matrix}/{Tensor}
  {Factorization} {Based} on {Low}-{Rank} {Approximation}.
\newblock IEEE Transactions on Signal Processing \textbf{60}(6), 2928--2940
  (2012).
\newblock \doi{10.1109/TSP.2012.2190410}.
\newblock \urlprefix\url{http://ieeexplore.ieee.org/document/6166354/}

\end{thebibliography}

\end{document}